\RequirePackage{fix-cm}
\documentclass[onecolumn]{svjour3questa}       
\usepackage{amsmath,amssymb}
\usepackage{mathptmx}      
\usepackage{hyperref}

\begin{document}

\title{Quasi-stationary distributions for queueing and other models}

\author{Phil Pollett}

\institute{Phil. Pollett \at
           School of Mathematics and Physics,
           The University of Queensland, Australia.
           \email{pkp@maths.uq.edu.au}
       }

\date{\today}

\maketitle

\newcommand{\Z}{\mathbb{Z}}

\section{Background}

We will examine questions concerning quasi-stationary behaviour in
evanescent processes. The idea has its origins in biological modelling,
where typically we are interested in limiting behaviour conditional on
non-extinction. For queueing processes we are typically interested in
the behaviour within a busy period in stable or near stable
queues~\cite{Kij93a,Kyp72b}, or, for unstable queues, prior to last exit
from the empty state~\cite{CHP00}.

Let $(X(t),\, t\geq 0)$ be a Markov chain in continuous time whose state
space $S=\{0\} \cup C$ consists of transient states~$C=\{1,2,\dots\}$
and an absorbing state~0. Let $Q=(q_{ij},\, i,j\in S)$ be the $q$-matrix of
transition rates, assumed to be stable, conservative and {\em regular\/}, 
so that there is a unique transition function $P(t)=(p_{ij}(t),\, i,j\in
S)$ associated with $Q$, and $p_{ij}(t)=\Pr(X(t)=j|X(0)=i)$. We
assume that $C$ is irreducible, and that $0$ is reached with
probability~$1$ from any state in $C$. Thus, in particular, for $i,j\in
C$, $p_{ij}(t)\to 0$ as $t\to\infty$, and, for $i\in C$, $p_{i0}(t) \to
1$. We are interested in the limit of the ratio
\begin{equation}
 \frac{p_{ij}(t)}{1-p_{i0}(t)} = \Pr(X(t)=j| X(t) \neq 0,\, X(0)=i),
 \qquad i,j\in C,
 \label{PKP2}
\end{equation}
called a {\em limiting conditional distribution\/} (LCD). The first
general results on LCDs~\cite{SV66} were facilitated by
a finer classification of transient states~\cite{Kin63a} and a common
rate at which the transition probabilities decay: there is a
$\lambda\geq 0$, called the {\em decay parameter\/} (of $C$), such that
$t^{-1} \log p_{ij}(t) \to -\lambda$ as $t\to \infty$, for all $i,j\in
C$. $C$ is then {\em $\lambda$-recurrent\/} or {\em
$\lambda$-transient\/} according to whether $\int_{0}^{\infty}
e^{\lambda t} p_{ij}(t)\, dt$ diverges or converges (for some, and then
all, $i,j\in C)$, and, when $C$ is $\lambda$-recurrent, it is {\em
$\lambda$-positive\/} or {\em $\lambda$-null\/} according to whether
the limit
$\lim_{t\to\infty} e^{\lambda t} p_{ij}(t)$ is positive or zero (for
some, and then all, $i,j\in C)$. 
When positive, its value is determined
by $(m_i,\, i\in C)$ and $(x_i,\, i\in C)$ satisfying

\begin{equation}
  \sum_{i\in C} m_i p_{ij}(t) = e^{-\lambda t} m_j
  \quad \text{and} \quad
  \sum_{j\in C} p_{ij}(t) x_j = e^{-\lambda t} x_i.
\label{PKP3}
\end{equation}
Unique positive solutions to~(\ref{PKP3}) 
(called the $\lambda$-invariant measure and vector, respectively)
are guaranteed when $C$ is $\lambda$-recurrent.
$C$ is then $\lambda$-positive if and only if 
$A^{-1}:=\sum_{k\in C} m_k x_k <\infty$, whence
$\lim_{t\to\infty} e^{\lambda t} p_{ij}(t) = A x_i m_j$.
So, one can see, at least formally from~(\ref{PKP2}), 
that, since $p_{i0}(t) = 1-\sum_{j\in C} p_{ij}(t)$,
$$
 \Pr(X(t)=j| X(t) \neq 0,\, X(0)=i) = 
 \frac{e^{\lambda t} p_{ij}(t)}{e^{\lambda t} \sum_{k\in C} p_{ik}(t)} 
 \to \frac{m_j}{\sum_{k\in C} m_k}.
$$
$\lambda$-positivity is indeed sufficient~\cite{Ver69}, the limit
taken to be~$0$ when $\sum_{k\in C} m_k=\infty$.

One might think this completes the picture. However,
$\lambda$-positivity is {\em not necessary\/} for the existence of an
LCD (see the example below). Further, the decay parameter cannot usually
be determined from $Q$, and $\lambda$-positivity cannot usually be
checked from~$Q$. Of course $p_{ij}(t)$ is seldom available explicitly, 
but there are ``$q$-matrix versions'' of~(\ref{PKP3}),
\begin{equation}
 \sum_{i\in C} m_i q_{ij}= -\lambda  m_j
  \quad \text{and} \quad
 \sum_{j\in C} q_{ij} x_j = -\lambda x_i,
\label{PKP4}
\end{equation}
and positive solutions to~(\ref{PKP4}) satisfy~(\ref{PKP3}) under
conditions that are easy to check~\cite{Pol86}.

\smallskip
\noindent
{\bf Example}\ \ 
Consider the M/M/1 queue with arrival rate~$p$ and departure rate~$q$
($>p$) modified so that it is killed when the queue size first reaches~$0$. 
Set $a = p+q$, $b=\sqrt{p/q}$ $(<1)$, and $\theta= 2\sqrt{pq}$.
Seneta~\cite{Sen66a} showed that, as $t\to\infty$,
$$
p_{ij}(t) 
 =i\, j b^{j-i} 
 \frac{2 e^{-(a-\theta)t} }{\theta\sqrt{2\pi \theta}} 
 \left( \frac{1}{t^{3/2}} + O \left( \frac{1}{t^{5/2}}\right) \right),
 \qquad i,j \in C,
$$
which implies that $\lambda = a - \theta = p+q - 2\sqrt{pq}$ is 
the decay parameter, and
$$
 \lim_{t\to \infty} t^{3/2} e^{\lambda t} p_{ij}(t) 
 =i\, j b^{j-i} \frac{2}{\theta\sqrt{2\pi \theta}},
\qquad i,j\in C.
$$
Notice that this limit is of the form $A x_i m_j$, where $m_j =
j\beta^{j}$ and $x_i = i\beta^{-i}$ specify the unique positive
solutions to~(\ref{PKP4}), and $A>0$. Seneta also showed that
$$
 \lim_{t\to\infty} t^{3/2} e^{\lambda t} \left(1-p_{i0}(t)\right)
 =\frac{ ib^{-i}}{\lambda \sqrt{2\pi \theta}},
\qquad i\in C.
$$
Notice also that this limit is of the form $B x_i$, where $B>0$.
So, the LCD exists:
$$
 \lim_{t\to\infty} \frac{p_{ij}(t)}{1-p_{i0}(t)}
 = \lim_{t\to\infty} \frac{t^{3/2}e^{\lambda t}
  p_{ij}(t)}{t^{3/2}e^{\lambda t} (1-p_{i0}(t))}
  =(1-b)^2 j b^{j-1},
\qquad i,j\in C.
$$
Yet, $C$ is $\lambda$-transient: $\int_{0}^{\infty} e^{\lambda t}
p_{ij}(t)\, dt = 2i/\theta<\infty$.

\section{Speculation}

One way to approach the question of whether LCDs exist for
$\lambda$-transient chains is to characterise the smallest $\kappa>1$
such that $t^{\kappa} e^{\lambda t} p_{ij}(t)$ has a strictly positive
limit for all $i,j\in C$. Such a characterisation is presently
unavailable. I conjecture (for $\lambda$-null and $\lambda$-transient
chains) that (i) {\em when\/} such a $\kappa$ exists, it is the same for
all $i,j\in C$, that (ii) the limit is always of the form $Ax_i m_j$,
where $(m_i,\, i\in C)$ and $(x_i,\, i\in C)$ satisfy~(\ref{PKP3})
perhaps with an inequality ($\leq$), and (iii) there is a $\kappa_0\leq
\kappa$ such that $t^{\kappa_0} e^{\lambda t} (1-p_{i0}(t)) \to B x_i$.

\section{Discussion}

The approach proposed here contrasts with work on {\em discrete-time\/}
chains, which highlight the many complications in the theory of
$R$-transient chains~\cite{McDF17,Kes95}, and more in line with work on
algebraic transience~\cite{MS14} and sub-geometric convergence for
ergodic discrete-time Markov chains~\cite{DMS07} (which has become
important in the analysis of Markov chain Monte Carlo methods). The
condition that $\kappa$ exists is equivalent to requiring that the
function $g_{ij}(t):=e^{\lambda t} p_{ij}(t)$ is {\em regularly
varying\/} with index $\alpha=\kappa^{-1}$ (the same index of regular
variation for all $i$ and~$j$), or, equivalently, that there is a
$\kappa>1$, the same for all $i,j\in C$, such that $t^{\kappa}
g_{ij}(t)$ is {\em slowly varying\/}~\cite{Sen76}. Conjectures (i) and
(ii) are true for $\lambda$-null recurrent chains (Lemma~1
of~\cite{Pol01}), and indeed~(iii) with $\kappa_0=\kappa$ under an
additional condition.
Critical to the argument is that the $\lambda$-subinvariant measures and
vectors ((\ref{PKP3}) with the inequality) are unique and
$\lambda$-invariant. This is not true in the $\lambda$-transient case.
One might hope to adapt Kingman's arguments based on inequalities
derived from the Chapman-Kolmorogov equations~\cite{Kin63a}, but I
cannot see how. In addition to the example detailed above, the
conjectures are supported by several other contrasting models: the M/M/1
queue with $p>q$, the random walk on $\Z$ in continuous
time~\cite{Kin63a}, the birth-death immigration process~\cite{And91},
quasi-birth-death processes~\cite{BBLPPT97}, and various branching
models~\cite{AH83}. Interestingly, the critical Markov branching process
provides an example for which $\kappa=2$ and $\kappa_0=1$.

\bibliographystyle{abbrv}
\bibliography{Pollett}

\end{document}